\documentclass[preprint,12pt]{elsarticle}

\usepackage{amssymb}
\usepackage{amsthm}
\usepackage{amsmath}
\usepackage{graphics}
\usepackage{graphicx}
\usepackage{tabularx}
\usepackage{appendix}
\usepackage{tikz}
\usepackage{xcolor}
\usepackage{booktabs}
\usepackage{bm}
\usepackage{pgfplots}
\usepackage{comment}
\usepackage[hidelinks]{hyperref}
\usepackage{subfig}
\usepackage{multirow}

\usepackage{xcolor}
\usepackage{soul}

\renewcommand{\tilde}[1]{\widetilde{#1}}

\usepackage{lineno}

\pgfplotsset{compat=newest, width = 6.5 cm}
\usepgfplotslibrary{colormaps}
\pgfplotsset{
/pgfplots/colormap={jet}{rgb255(0cm)=(0,0,128) rgb255(1cm)=(0,0,255)
rgb255(3cm)=(0,255,255) rgb255(5cm)=(255,255,0) rgb255(7cm)=(255,0,0) rgb255(8cm)=(128,0,0)}
}

\begin{document}
\begin{frontmatter}

\title{Consistent lattice Boltzmann methods for the volume averaged Navier\textendash Stokes equations}

\author[inst0,inst1]{Fedor Bukreev\corref{cor1}}
\cortext[cor1]{Corresponding author: fedor.bukreev@kit.edu}

\author[inst0,inst2]{Stephan Simonis}
\author[inst0,inst1,inst2]{Adrian Kummerländer}
\author[inst0,inst1]{Julius Jeßberger}
\author[inst0,inst1,inst2]{Mathias J. Krause}

\address[inst0]{Lattice Boltzmann Research Group (LBRG),}
\address[inst1]{Institute of Mechanical Process Engineering and Mechanics (MVM),}
\address[inst2]{Institute for Applied and Numerical Mathematics (IANM), \\ Karlsruhe Institute of Technology (KIT), 76131 Karlsruhe, Germany}

\begin{abstract}
We derive a novel lattice Boltzmann scheme, which uses a pressure correction forcing term for approximating the volume averaged Navier\textendash Stokes equations (VANSE) in up to three dimensions. 
With a new definition of the zeroth moment of the Lattice Boltzmann equation, spatially and temporally varying local volume fractions are taken into account. 
A Chapman\textendash Enskog analysis, respecting the variations in local volume, formally proves the consistency towards the VANSE limit up to higher order terms. 
The numerical validation of the scheme via steady state and non-stationary examples approves the second order convergence with respect to velocity and pressure. 
The here proposed lattice Boltzmann method is the first to correctly recover the pressure with second order for space-time varying volume fractions.
\end{abstract}

\begin{keyword}
volume averaged Navier\textendash Stokes \sep lattice Boltzmann method \sep consistency \sep Chapman\textendash Enskog analysis
\end{keyword}

\end{frontmatter}

\section{Introduction}
\label{sec:introduction}

\noindent Multiphase simulations gain high demand in both industry and science, especially in the field of process technology, where complex reactions and phase transitions need to be calculated. 
Exemplary applications are liquid-solid, gas-solid and gas-liquid reactors \cite{reschetilowski}, phase separators and transport units \cite{sattler, bohnet}. 
Today, computational fluid dynamics (CFD) eases the design and optimization of such processes. 
Three different strategies for the simulation of complex flows akin to the above have been established, namely volume of fluid (VOF) methods, discrete element methods (DEM), and Eulerian multiphase methods \cite{hiltunen}. 
Whereas VOF methods track the phase interface in detail, the DEM approach calculates paths of discrete particles.

The Eulerian multiphase simulation schemes typically consider each phase as continuous and solve mass and momentum conservation equations for each of them. 
The coupling between phases is realized via volume averaging of the phase flow variables. 
The resulting volume averaged Navier\textendash Stokes equations (VANSE) \cite{GIDASPOW1994337} involve also the phase interaction forces in the momentum equation. 
Besides multiphase flows, the VANSE are suitable for the modeling of porous flows \cite{zhu}. 
In comparison to VOF, the Eulerian methods require lower computational resources in general. 
Similarly, the latter outperform the DEM if the amount of particles reaches billions or more. 

The Navier\textendash Stokes equations (NSE) and VANSE can be solved numerically in the discretized form with the finite difference method (FDM) \cite{PEPIOT2012104}, finite element method (FEM) or finite volume method (FVM) \cite{moukalled} on the macroscopic level or with the lattice Boltzmann methods (LBM) \cite{kruger}, which are based on mesoscopic kinetic theory \cite{hanel2004molekulare}. 
In LBM the fluid is considered as a quantity of colliding and streaming particles. 
The state of particles is described by a discretized particle distribution function (population) \textemdash the probability of the particle to be located at the regarded coordinates in the phase space. 
The equilibrium population is the Maxwellian distribution based on the equation of state. 
The collision and streaming of populations is described by a simplified version of the Boltzmann equation.  
Taking moments of one lattice cell leads to the macroscopic quantities density, velocity and pressure, respectively. 
Through the Chapman\textendash Enskog (CE) expansion \cite{li} or limit consistency \cite{simonis2022limit}, the lattice Boltzmann equation can be linked to the NSE. 
The most prolific feature of LBM is the suitability for parallelization due to explicitly local calculation of populations. 
Meanwhile, LBM has been found to provide advanced capabilities for the parallel simulation of turbulent flows \cite{simonis2022temporal,simonis2021linear,haussmann2019direct}, advection\textendash diffusion transport \cite{simonis2020relaxation,dapelo2021lattice-boltzmann}, and more specific photobioreactors \cite{mink2021comprehensive}, Flettner rotors \cite{simonis2022forschungsnahe} or Coriolis mass flow meters \cite{haussmann2021fluid-structure}.
As a paragon of effectiveness of the LBM, the comparison between the open-source software packages OpenLB \cite{krause,kummerlander2022olb15} and OpenFOAM shows 32 times faster computation time of the former by the in-cylinder flow test \cite{krause, tuprints}. 

Particular LBM for the solution of VANSE were developed by several authors. 
Ansatz of Guo et al. \cite{guo} for flows through porous media is a discretization of the Darcy\textendash Lapwood\textendash Brinkman equation. 
Unfortunately, this realization is only valid for temporally and spatially constant void fractions. 
Blais et al.~\cite{blais} proposed a scheme, which is based on the method of moments, where first the population moments necessary for the VANSE are chosen and after that the equilibrium distribution is composed. 
The volume fraction is implemented only into the zeroth population, what makes the pressure calculation more stable, but allows application of this model only by volume fractions above $0.5$. 
This model fits the majority of porous flows, but is not universally applicable for all multiphase flows. 
Although Höcker et al.~\cite{hocker} and Maier et al.~\cite{Maier2021_1000132643} correct the zeroth moment on the lattice Boltzmann level, the CE expansion of this method in case of strongly varying local volume fractions is not fulfilled. 
The simplest and most uniform VANSE LBM is suggested by Zhang et al.~\cite{zhang}. 
The method fits cases with temporally and spatially varying volume fractions except for the pressure distribution. 
To the knowledge of the present authors, the pressure discrepancy in \cite{zhang} due to an inconsistent zeroth moment interpretation in the there performed CE expansion. 
This in turn leads to a density calculation which changes pressure correction forces and the pressure itself. 
Based on the preceding approaches, the present work proposes a consistent way of the numerical VANSE solution with lattice Boltzmann methods for one, two, and three dimensions.

The paper is structured as follows. 
First, the principles of VANSE and the corresponding LBM scheme are derived in Section \ref{sec:meth}. 
In particular, the novel population moments are presented and locally varying void fractions are taken into account. 
In Section \ref{sec:numerics}, the validation of the new correction is performed on stationary and transient examples with spatially changing volume fractions between $0.1$ and $0.9$. 
The numerical results suggest a second order convergence of flow velocity and pressure. 
Section \ref{sec:conc} draws conclusions and suggests future research. 
At last, the CE expansion, formally proving the approximation of the VANSE with the present LBM up to higher order terms, is detailed in \ref{sec:appendix}.

\section{Methodology}\label{sec:meth}

\subsection{Volume averaged Navier\textendash Stokes equations}
\noindent If subgrid particles are contained in the regarded control volume, any quantity of a fluid phase can be adjusted to the whole volume, which also includes these particles. 
Below, this adjustment is called volume averaging, denoted with $\langle \cdot \rangle$, and can be written as follows for any fluid quantity $q^{\flat}$, where $\flat$ indicates the corresponding phase. 

Let \(V\) denote the overall volume and \(V^{\flat}\) the volume which is occupied by phase $\flat$, hence 
\begin{linenomath}\begin{align}
    \sum\limits_{\flat} V^\flat = V.
\end{align}\end{linenomath}
The ratio of these volumes is defined via the respective void fraction 
\begin{linenomath}\begin{align}
    \phi^{\flat} &= \frac{V^{\flat}}{V}.
\end{align}\end{linenomath}
In the following, volume averaged scalars are denoted with $\tilde{\cdot}$ and volume averaged vectors with $\overline{\cdot}$.
Thus, for scalars \(q\) and vectors \(\bm{s}\) we define
\begin{linenomath}\begin{align}
    \phi^{\flat} \tilde{q^{\flat}} &= \langle q^{\flat} \rangle \equiv \frac{1}{V} \int\limits_{V^{\flat}} q^{\flat} \,\mathrm{d}V, \\
    \phi^{\flat} \tilde{\rho^{\flat}} \overline{\bm{s}^{\flat}} &= \langle \rho^{\flat} \bm{s}^{\flat} \rangle.
\end{align}\end{linenomath}
By volume averaging all terms of the NSE, the VANSE are deduced \cite{hiltunen, GIDASPOW1994337}
\begin{linenomath}\begin{align}
    \partial_t (\phi^{\flat} \tilde{\rho^{\flat}})+ \bm{\nabla} \cdot (\phi^{\flat} \tilde{\rho^{\flat}} \overline{\bm{u}^{\flat}}) & = 0, \label{eq:vansEquMass} \\    
    \partial_t (\phi^{\flat} \tilde{\rho^{\flat}} \overline{\bm{u}^{\flat}})+ \bm{\nabla} \cdot (\phi^{\flat} \tilde{\rho^{\flat}} \overline{\bm{u}^{\flat}} \overline{\bm{u}^{\flat}})+ \phi^{\flat} \bm{\nabla} \tilde{p} & = \nu \bm{\nabla} \cdot (\phi^{\flat} \tilde{\rho^{\flat}} (\bm{\nabla} \overline{\bm{u}^{\flat}} + \overline{\bm{u}^{\flat}} \bm{\nabla})) + \phi^{\flat} \overline{\bm{F}^{\flat}}, \label{eq:vansEquMom}
\end{align}\end{linenomath}
where $\tilde{\rho^{\flat}}$ and $\overline{\bm{u}^{\flat}}$ denote the volume averaged versions of the fluid density and the velocity, respectively. 
The pressure $\tilde{p}$ is common for all phases in the system.

\subsection{Lattice Boltzmann scheme for volume averaged Navier\textendash Stokes equations}

\begin{figure}[h]
  \centering
  \includegraphics[scale=1]{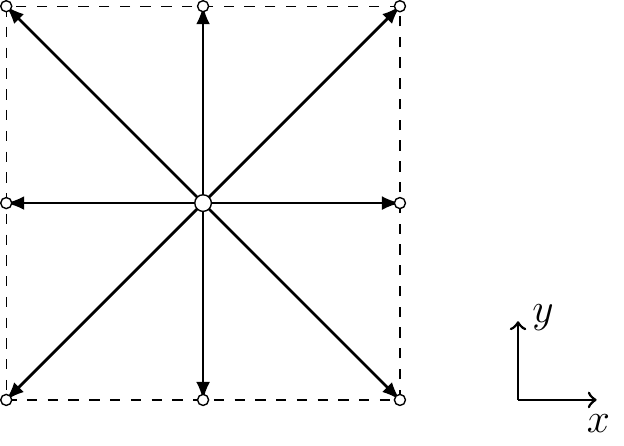}
 \caption{Schematic view of the \(D2Q9\) discrete velocity set.}
 \label{fig:D2Q9_lattice}
\end{figure}

\begin{figure}
\centering
  \includegraphics[scale=1]{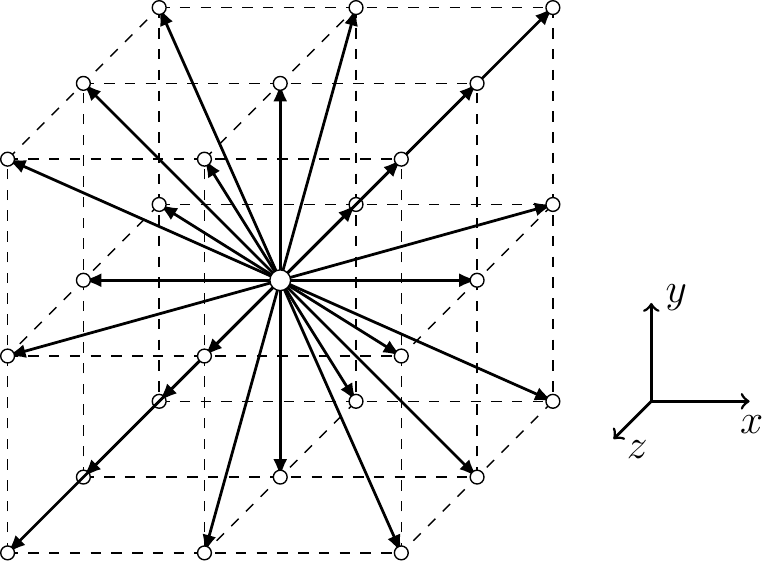}
  \caption{Schematic view of the \(D3Q27\) discrete velocity set.}
  \label{fig:D3Q27_lattice}
\end{figure}

\noindent In the following, equations \eqref{eq:vansEquMass} and \eqref{eq:vansEquMom} are approximated with an LBM based on Bhatnagar\textendash Gross\textendash Krook (BGK) collision \cite{bgk} and Guo et al. forcing \cite{guozhaoli} on two- and three-dimensional \(D2Q9\) and \(D3Q27\) lattices. 
One-dimensional stencils are also discussed but not focused here. 
The discrete velocity sets are visualized in Figures \ref{fig:D2Q9_lattice} and \ref{fig:D3Q27_lattice}. 
The corresponding discretization parameters are given in Tables \ref{tab:latt_par} and \ref{tab:latt_3d}, respectively.

\begin{table}[h]
\centering
\begin{tabular}{ccc}
    \toprule
    Directions \(i\) & Normalized lattice velocity $\bm{\xi}_{i}$ & Lattice weights $w_{i}$ \\
    \midrule
    $0$			& $(0, 0)$		            & $4/9$ \\
    $1,2,3,4$		& $(\pm 1, 0), (0, \pm 1)$	& $1/9$ \\
    $5,6,7,8$		& $(\pm 1, \pm 1)$		    & $1/36$ \\
    \bottomrule
\end{tabular}
\caption{Lattice discretization parameters of \(D2Q9\).}
    \label{tab:latt_par}
\end{table}

\begin{table}[h]
\centering
\begin{tabular}{ccc}
    \toprule
    Directions \(i\) & Normalized lattice velocity $\bm{\xi}_{i}$ & Lattice weights $w_{i}$ \\
    \midrule
    $0$			& $(0, 0)$		            & $8/27$ \\
    $1,2,\ldots,6$		& $(\pm 1, 0, 0), (0, \pm 1, 0), (0, 0, \pm 1)$	& $2/27$ \\
    $7,8,\ldots,18$		& $(\pm 1, \pm 1, 0), (\pm 1, 0, \pm 1), (0, \pm 1, \pm 1)$		    & $1/54$ \\
    $19,20,\ldots,26$			& $(\pm 1, \pm 1, \pm 1)$		            & $1/216$ \\
    \bottomrule
\end{tabular}
\caption{Lattice discretization parameters of \(D3Q27\).}
    \label{tab:latt_3d}
\end{table}

Unless stated otherwise, \(i=0,1,\ldots,26\) denotes the population index. 
The space-time discrete lattice Boltzmann equation (LBE) reads
\begin{linenomath}\begin{align}\label{eq:lbe}
    f_{i}(\bm{x}+\bm{\xi}_{i} \triangle t, t+\triangle t) = f_{i}(\bm{x},t) + \frac{\triangle t}{\tau}(f_{i}^{eq}(\bm{x},t) - f_{i}(\bm{x},t)) + \Omega^{F}_{i} .
\end{align}\end{linenomath}
The equilibrium particle distribution function used by Zhang et al. \cite{zhang} as well as by Höcker et al. \cite{hocker} and Maier \cite{Maier2021_1000132643} is simple, universal for all populations from 0 to 26, and stable for all possible volume fraction values. 
It is the common, second order truncated Maxwell equilibrium, multiplied with the local volume fraction.
\begin{linenomath}\begin{align}\label{eq:equilibrium}
 f_{i}^{\mathrm{eq}}(\bm{x},t) = w_{i} \tilde{\rho^{\flat}} \phi^{\flat} \Bigl(1 + \frac{\xi_{i\alpha} \overline{u_{\alpha}^{\flat}}}{c_{s}^{2}} + \frac{(\xi_{i \alpha} \xi_{i \beta} - c_{s}^{2} \delta_{\alpha \beta}) \overline{u_{\alpha}^{\flat}} \overline{u_{\beta}^{\flat}}}{2c_{s}^{4}} \Bigr)
\end{align}\end{linenomath}
After the first time step, $\tilde{\rho^{\flat}} \phi^{\flat}$ is replaced by the zeroth population moment $\sum_{i} f_{i}$.
The standard LBM presupposes the constant density of the fluid, which is typically fulfilled e.g. in multiphase or porous flows. 
In contrast, if the constant density is multiplied with the spatially and temporally varying volume fraction, the result is not constant anymore.  
The density in lattice units takes usually the value of 1, whereas the volume fraction can vary between 0 and 1, such that the effective density considered here in turn is varying also between 0 and 1. 
Taking into account the streaming of effective densities along the lattice directions, the new form of the equilibrium distribution function is then after the first collision
\begin{linenomath}\begin{align}
 f_{i}^{\mathrm{eq}}(\bm{x},t) = w_{i} \tilde{\rho^{\flat}} \Bigl( \int_V \phi^{\flat}(\bm{x},t) \mathrm{d}V \Bigr) \Bigl(1 + \frac{\xi_{i\alpha} \overline{u_{\alpha}^{\flat}}}{c_{s}^{2}} + \frac{(\xi_{i \alpha} \xi_{i \beta} - c_{s}^{2} \delta_{\alpha \beta}) \overline{u_{\alpha}^{\flat}} \overline{u_{\beta}^{\flat}}}{2c_{s}^{4}} \Bigr) . \label{eq:feq}
\end{align}\end{linenomath}

Based on that, we define the effective density and velocity as 
\begin{linenomath}\begin{align}
 \tilde{\rho^{\flat}} & = \frac{\sum_{i} f_{i}}{\int_V \phi^{\flat}(\bm{x},t) \mathrm{d}V}, \label{eq:density} \\
 \overline{\bm{u}^{\flat}} & = \frac{\sum_{i} \bm{\xi}_{i} f_{i}}{\sum_{i} f_{i}} + \frac{\triangle t}{2} \frac{\sum_{k} \overline{\bm{F}_{k}}}{\sum_{i} f_{i}}, \label{eq:velocity}
\end{align}\end{linenomath}
respectively. 
The density definition uses for the volume fraction integration the neighboring cell data, which is considered further below. 
Due to Guo et al. forcing scheme \cite{guozhaoli}, the velocity $\overline{\bm{u}^{\flat}}$ contains the sum of forces used in the example $\sum_{k} \overline{\bm{F}_{k}}$. 
Further, the forcing term is defined as
\begin{linenomath}\begin{align}
 \Omega_{i}^{F} = \Bigl(1 - \frac{\triangle t}{2\tau} \Bigr) w_{i} \Bigl(\frac{\xi_{i\alpha}}{c_{s}^{2}} + \frac{(\xi_{i \alpha} \xi_{i \beta} - c_{s}^{2} \delta_{\alpha \beta}) \overline{u_{\beta}^{\flat}}}{c_{s}^{4}} \Bigr)  \sum\limits_{k} \overline{F_{k\alpha}}. 
\end{align}\end{linenomath}
 The sum \(\sum_{k} \overline{F_{k\alpha}} \) includes the phase interaction forces and the pressure correction force proposed by Zhang et al. \cite{zhang}
\begin{linenomath}\begin{align}\label{eq:pressForce}
 \overline{\bm{F}_{\mathrm{PC}}} = \tilde{p} \bm{\nabla} \phi^{\flat} = \tilde{\rho^{\flat}} c_{s}^{2} \bm{\nabla} \phi^{\flat}.
\end{align}\end{linenomath}
This correction force adjusts the pressure term in the momentum equation, which is $\bm{\nabla} (\phi^{\flat} \tilde{p})$ according to the CE expansion of Zhang et al. equilibrium particle distribution and should be $\phi^{\flat} \bm{\nabla} \tilde{p}$ as in VANSE. 
The phase interaction forces are for example in the case of a particle-laden flow given by the drag, lift, gravity, virtual mass and turbulence interaction forces. 
These interaction forces are not considered in the present work due to the focus on model validation. 
Note that the consistent incorporation of the neglected forces can be done with Guo et al. forcing scheme alongside the pressure correction. 
Hence, without loss of generality we assume that \(\sum_{k} \overline{\bm{F}_{k}} = \overline{\bm{F}_{\mathrm{PC}}}\). 
Further, the gradient of volume fraction appearing in \eqref{eq:pressForce} is discretized through central differences, thus for example in two dimensions
\begin{linenomath}\begin{align}
 \bm{\nabla} \phi^{\flat} \approx \frac{1}{2\triangle x} \begin{pmatrix} \phi_{x+1}^{\flat} - \phi_{x-1}^{\flat} \\ \phi_{y+1}^{\flat} - \phi_{y-1}^{\flat} \end{pmatrix}. 
\end{align}\end{linenomath}

The above-mentioned effective density is part of the equilibrium distribution function, and hence propagates from and to the neighbor lattice cells (cf. \eqref{eq:equilibrium} $\rightarrow$ cf. \eqref{eq:feq}), such that volume fraction becomes integrated over the cell volume \(\tilde{\rho^{\flat}} \phi^{\flat} \rightarrow \tilde{\rho^{\flat}}\int_V \phi^{\flat}(\bm{x},t) \mathrm{d}V\). 
Each cell contains own distinct effective density and different density values at the interfaces, calculated by integration with the neighbor cells effective densities. 
For the discretized integral calculations we use quadrature rules 
\begin{linenomath}\begin{align}\label{eq:quad}
 \int_V \phi^{\flat}(\bm{x},t) \mathrm{d}V &= \sum_i^{N} \varpi_i(N) \phi^{\flat}(\bm{x} - \bm{\xi}_{i} \triangle t,t)
\end{align}\end{linenomath}
which are rearranged to 
\begin{linenomath}\begin{align}
\sum_i^{N} \varpi_i(N) \phi^{\flat}(\bm{x} - \bm{\xi}_{i} \triangle t,t) & = \left( \varpi_{i \neq 0}(N) \bm{\nabla}^{2} \phi^{\flat} + \phi^{\flat} \right) \left( \bm{x}, t\right), \label{eq:latDiffRuleIIre} 
\end{align}\end{linenomath}
respectively. 
The number of quadrature points $N$ is dependent on the void fraction variation directions number. 
Hereby the diagonal directions are not considered. 
In particular, the volume fraction integration is performed on the $D1Q3$ lattice if the volume fraction changes only in one direction, on the $D2Q5$ lattice if in two and on $D3Q7$ if in all three directions. 
In \eqref{eq:quad} and \eqref{eq:latDiffRuleIIre} $N$ is equal to $Q$. 
It is to be noted that the weighting factors, which are listed in Table \ref{tab:weights}, do not conform to the weights of a discrete velocity set.

\begin{table}[h]
\centering
\begin{tabular}{ccc}
    \toprule
    Dimensions      & $\varpi_0$            & $\varpi_{i \neq 0}$ \\
    \midrule
    $d=1$			& $1/2$		            & $1/4$ \\
    $d=2$		    & $1/3$	                & $1/6$ \\
    $d=3$		    & $1/6$		            & $5/36$ \\
    \bottomrule
\end{tabular}
\caption{Quadrature weights for void fraction integration over a lattice cell.}
    \label{tab:weights}
\end{table}

The equilibrium moments with varying local volume fractions are thus computed via \eqref{eq:feq}, \eqref{eq:density}, \eqref{eq:velocity}, and \eqref{eq:latDiffRuleIIre} in a separately regarded lattice cell in the pre-collision state to 
\begin{linenomath}\begin{align}
 M^{\mathrm{eq}}_{0} &= \sum_{i} f_{i}^{\mathrm{eq}} = \tilde{\rho^{\flat}} \Bigl(\phi^{\flat} (\bm{x},t) + \varpi_{i \neq 0}(N) \bm{\nabla}^{2} \phi^{\flat} \Bigr), \label{eq:M0} \\
 M_{1\alpha}^{\mathrm{eq}} &= \sum_{i} \xi_{i\alpha} f_{i}^{\mathrm{eq}}  = \tilde{\rho^{\flat}} \Bigl(\phi^{\flat} (\bm{x},t) + \varpi_{i \neq 0}(N) \bm{\nabla}^{2} \phi^{\flat} \Bigr) \overline{\bm{u}^{\flat}}, \\
 M_{2\alpha\beta}^{\mathrm{eq}} &= \sum_{i} \xi_{i \alpha} \xi_{i \beta} f_{i}^{\mathrm{eq}} = \tilde{\rho^{\flat}} \Bigl(\phi^{\flat} (\bm{x},t) + \varpi_{i \neq 0}(N) \bm{\nabla}^{2} \phi^{\flat} \Bigr) \overline{u_{\alpha}^{\flat}} \overline{u_{\beta}^{\flat}} \nonumber \\
 &+  \tilde{\rho^{\flat}} c_s^{2} \Bigl(\phi^{\flat} (\bm{x},t) + \varpi_{i \neq 0}(N) \bm{\nabla}^{2} \phi^{\flat} \Bigr), \\ 
 M_{3\alpha\beta\gamma}^{\mathrm{eq}} &= \sum_{i} \xi_{i \alpha} \xi_{i \beta} \xi_{j \gamma} f_{i}^{\mathrm{eq}} = \tilde{\rho^{\flat}} \Bigl(\phi^{\flat} (\bm{x},t) + \varpi_{i \neq 0}(N) \bm{\nabla}^{2} \phi^{\flat} \Bigr) \overline{\bm{u}^{\flat}} \delta_{\alpha \beta \gamma}. \label{eq:M3}
\end{align}\end{linenomath}
Finally, using these moments, a CE expansion (see \ref{sec:appendix}) of the above proposed lattice Boltzmann scheme yields formal consistency towards the VANSE \eqref{eq:vansEquMass}, \eqref{eq:vansEquMom}.

\section{Numerical validation}\label{sec:numerics}
\noindent The numerical validation of the proposed LBM for VANSE is performed on a stationary and a transient example. 
Both examples are built with the method of manufactured solutions (MMS) \cite{roache}. 
Thereby, analytical functions for volume fraction, fluid velocity and pressure are chosen, s.t. they fulfill the mass conservation law of VANSE. 
For these fixed functions $\phi^{\flat},\overline{\bm{u}^{\flat}}, \tilde{p_{i}}$, the MMS force is calculated with central finite differences to 
\begin{linenomath}\begin{align}
    \bm{F}_{\mathrm{MMS}} = \partial_t \left(\phi^{\flat} \tilde{\rho^{\flat}} \overline{\bm{u}^{\flat}}\right) & + \bm{\nabla} \cdot \left(\phi^{\flat} \tilde{\rho^{\flat}} \overline{\bm{u}^{\flat}} \overline{\bm{u}^{\flat}}\right)+ \phi^{\flat} \bm{\nabla} \tilde{p_{i}} \nonumber \\
    & - \nu \bm{\nabla} \cdot \left(\phi^{\flat} \tilde{\rho^{\flat}} \left(\bm{\nabla} \overline{\bm{u}^{\flat}} + \overline{\bm{u}^{\flat}} \bm{\nabla}\right)\right) ,
\end{align}\end{linenomath}
 including all terms of the momentum equation.
This force is used as forcing term in the LBE \eqref{eq:lbe} together with the pressure correction force
\begin{linenomath}\begin{align}
    \sum_{k} \overline{\bm{F}_{k}} = \overline{\bm{F}_{\mathrm{MMS}}} + \overline{\bm{F}_{\mathrm{PC}}}.
\end{align}\end{linenomath}
The examples are evaluated through several error measurements. 
The errors correspond to $L^1$-, $L^2$- and $L^{\infty}$-norms over nodal values of velocity and pressure deviations between the simulated and the prescribed data \cite{oberkampf}, i.e.
\begin{linenomath}\begin{align}
    r_{L^{1}} \left( q^{\flat} \right) &= \frac{1}{N_{\mathrm{node}}} \sum_{c=1}^{N_{\mathrm{node}}} \left\vert q_{c}^{\flat} - q_{c}^{\flat,\star} \right\vert , \\
    r_{L^{2}} \left( q^{\flat} \right) &= \sqrt{ \frac{1}{N_{\mathrm{node}}} \sum_{c=1}^{N_{\mathrm{node}}} \left\vert q_{c}^{\flat} - q_{c}^{\flat,\star} \right\vert^{2}}, \\
    r_{L^{\infty}} \left( q^{\flat} \right) &= \max\limits_{c = 1,..,N_{\mathrm{node}}} \left\vert q_{c}^{\flat} - q_{c}^{\flat,\star} \right\vert,
\end{align}\end{linenomath}
respectively, where \(q^{\flat,\star}\) denotes the corresponding analytical solution. 

\begin{figure}
\centering
\begin{tabular}{cc}
    \multirow{2}{*}[9.5em]{
    \includegraphics[width=0.6\linewidth]{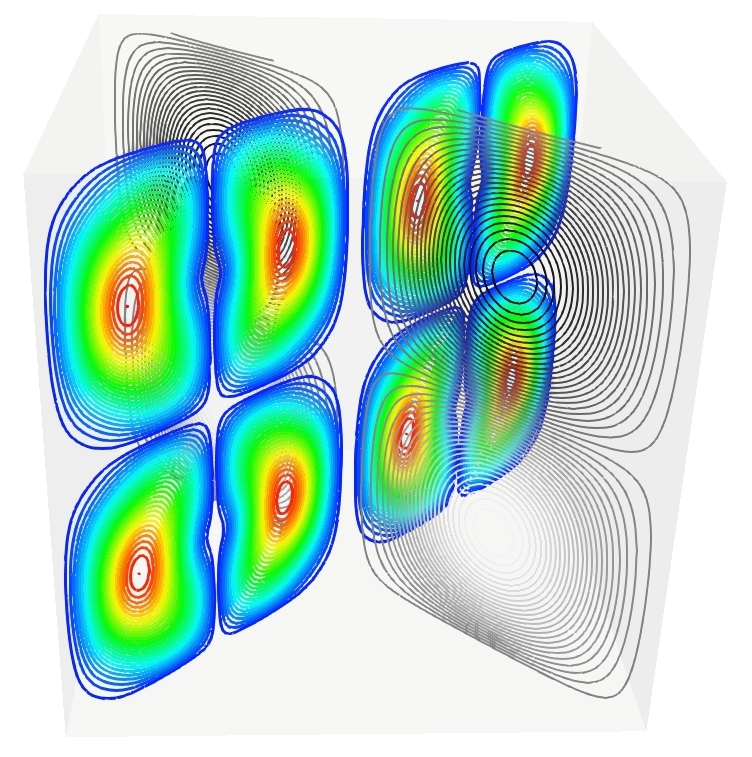}
	} 
    &
	\includegraphics[scale=1]{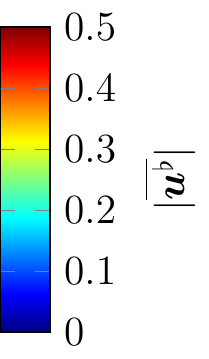} \\ 
    &
	\includegraphics[scale=1]{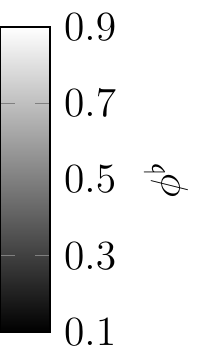}
\end{tabular}    
    \caption{Stationary three-dimensional velocity and porosity distribution of Example \ref{ex:stat3d}.}
    \label{fig:3D_stationary}
\end{figure}

\begin{figure}
\centering
\begin{tabular}{cc}
    \multirow{2}{*}[9.5em]{
    \includegraphics[width=0.6\linewidth]{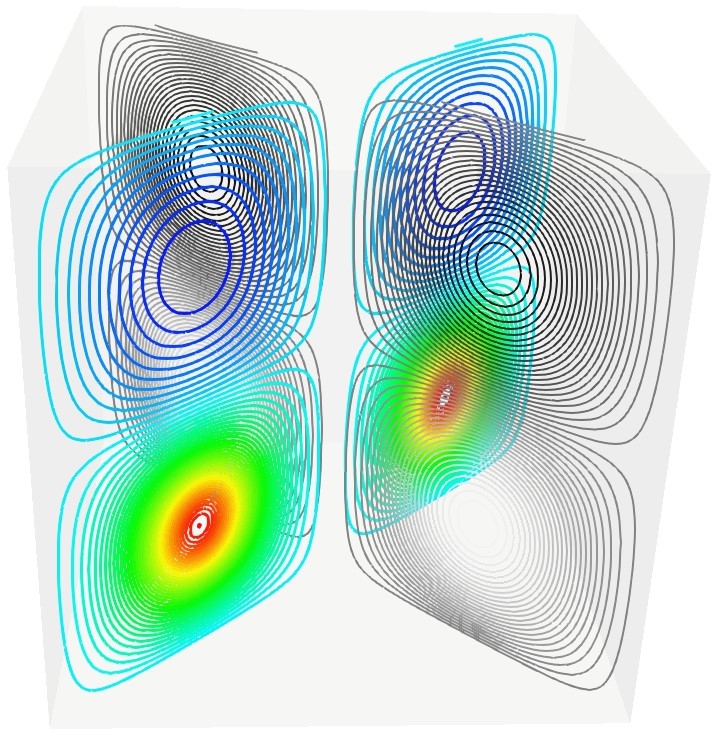}
	} 
    &
    \includegraphics[scale=1]{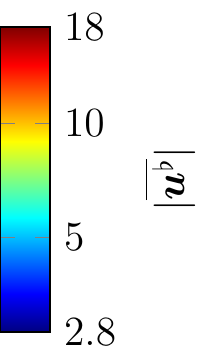} \\ 
    &
	\includegraphics[scale=1]{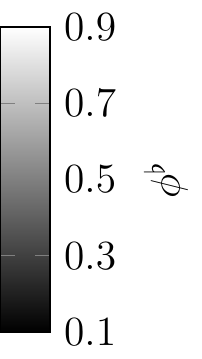}
\end{tabular}    
    \caption{Transient three-dimensional velocity and porosity distribution of Example \ref{ex:tran3d}.}
    \label{fig:3D_transient}
\end{figure}

\begin{figure}[ht!]
\centerline{
\subfloat[Velocity error]{
	\includegraphics[scale=1]{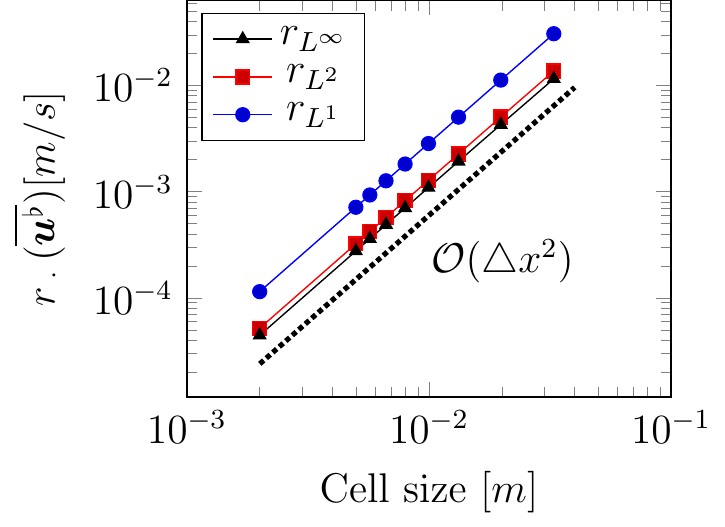}}
\subfloat[Pressure error]{
	\includegraphics[scale=1]{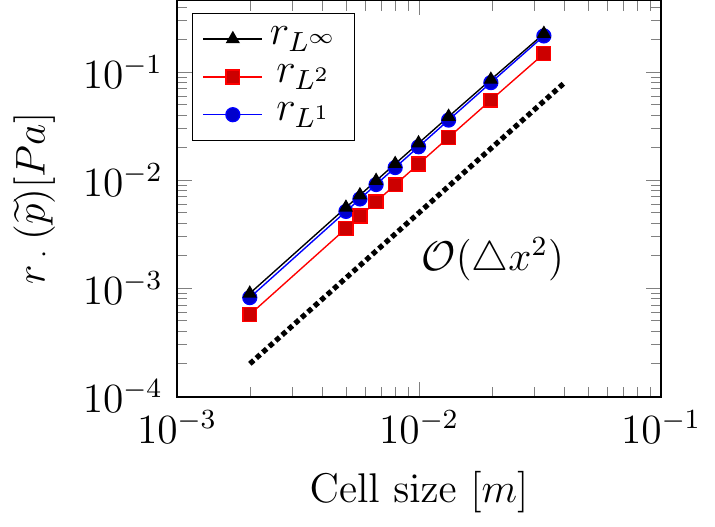}}
}
\caption{Error measurements for (a) velocity and (b) pressure of the stationary two-dimensional Example \ref{ex:stat2d}.}
\label{fig:stationary_plots_2d}
\end{figure}

\begin{figure}[ht!]
\centerline{
\subfloat[Velocity error]{
	\includegraphics[scale=1]{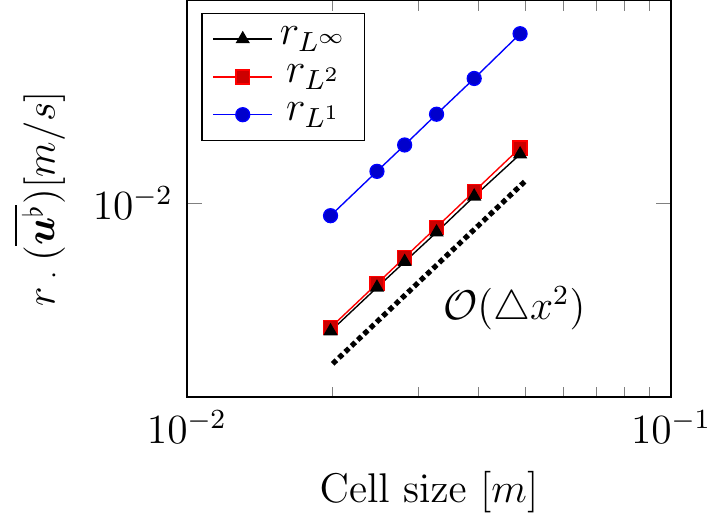}}
\subfloat[Pressure error]{
	\includegraphics[scale=1]{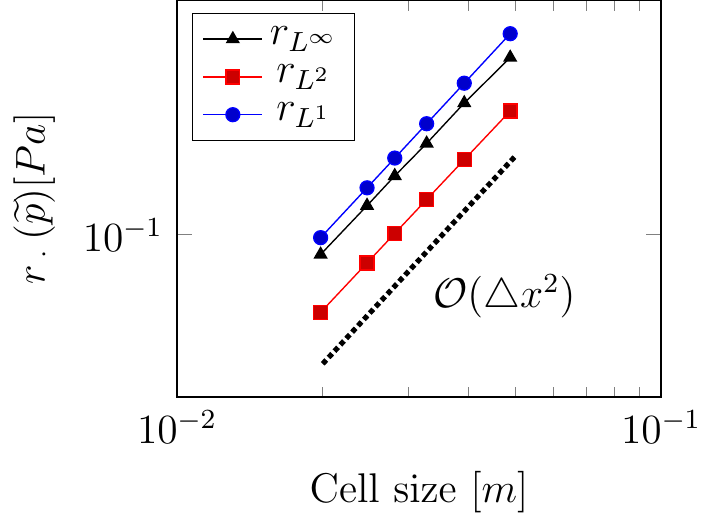}}
}
\caption{Error measurements for (a) velocity and (b) pressure of the stationary three-dimensional Example \ref{ex:stat3d}.}
\label{fig:stationary_plots_3d}
\end{figure}

\begin{figure}[ht!]
\centerline{
\subfloat[Velocity error]{
	\includegraphics[scale=1]{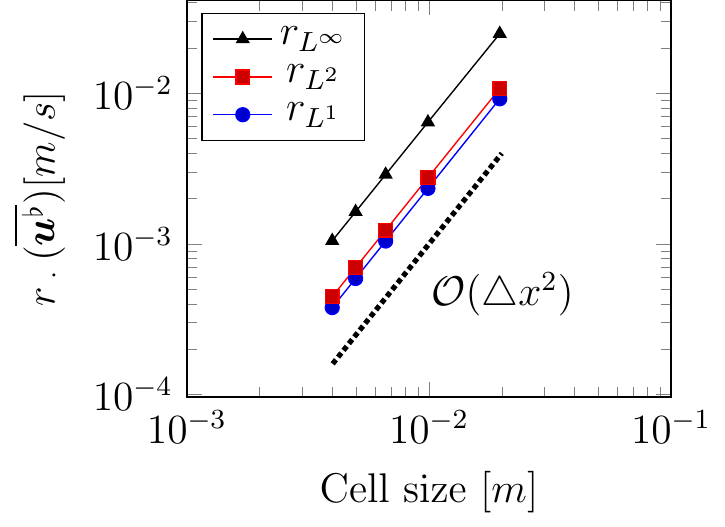}}
\subfloat[Pressure error]{
	\includegraphics[scale=1]{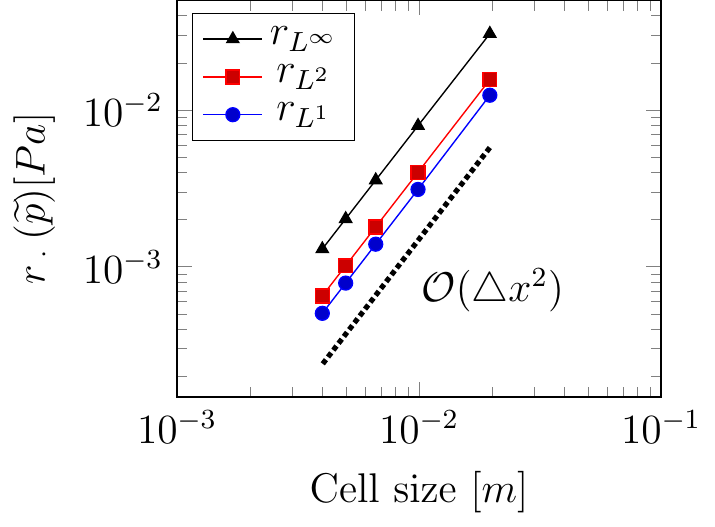}}
}
\caption{Error measurements for (a) velocity and (b) pressure of the transient one-dimensional Example \ref{ex:tran1d}.}
\label{fig:transient_plots_1d}
\end{figure}

\begin{figure}[ht!]
\centerline{
\subfloat[Velocity error]{
	\includegraphics[scale=1]{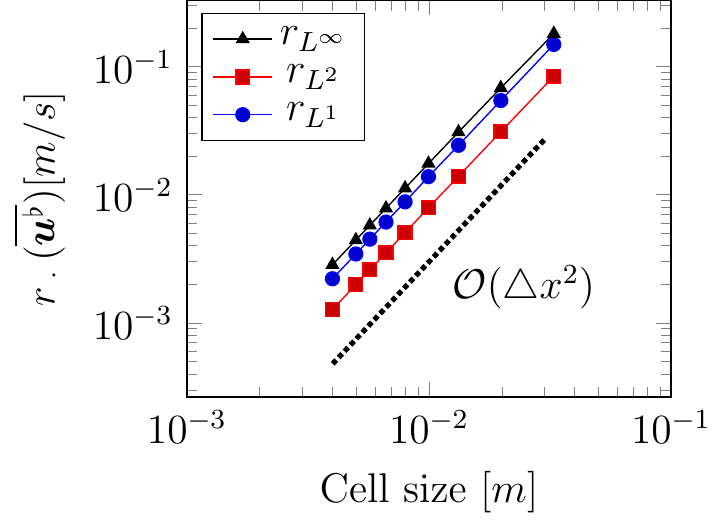}}
\subfloat[Pressure error]{
	\includegraphics[scale=1]{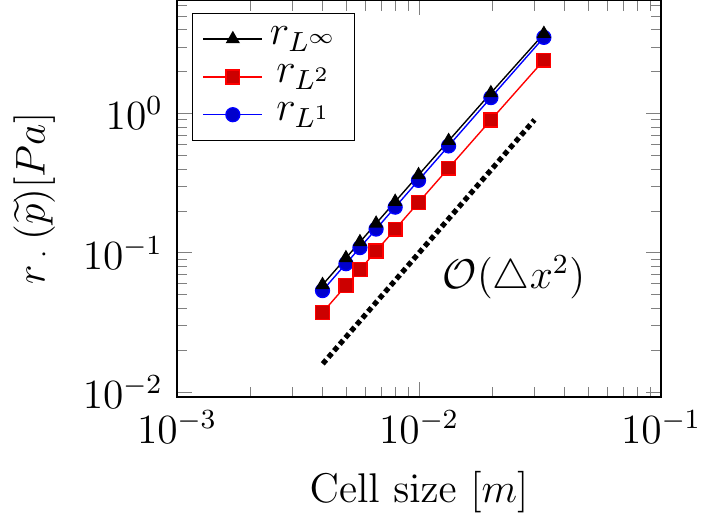}}
}
\caption{Error measurements for (a) velocity and (b) pressure of the transient two-dimensional Example \ref{ex:tran2d}.}
\label{fig:transient_plots_2d}
\end{figure}

\begin{figure}[ht!]
\centerline{
\subfloat[Velocity error]{
	\includegraphics[scale=1]{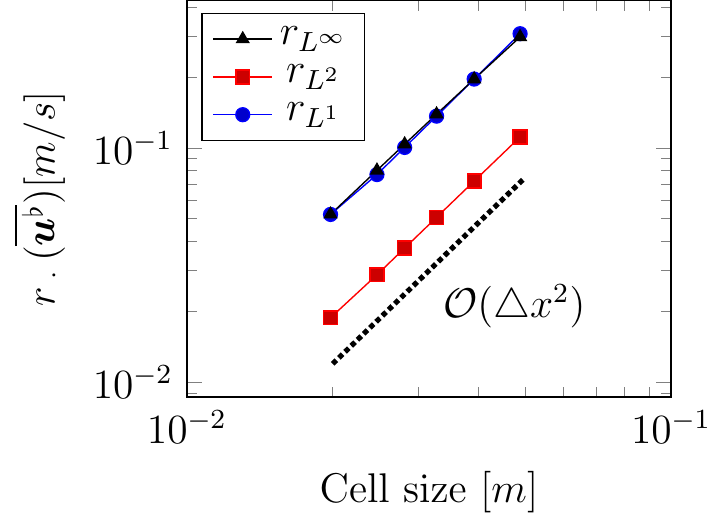}}
\subfloat[Pressure error]{
	\includegraphics[scale=1]{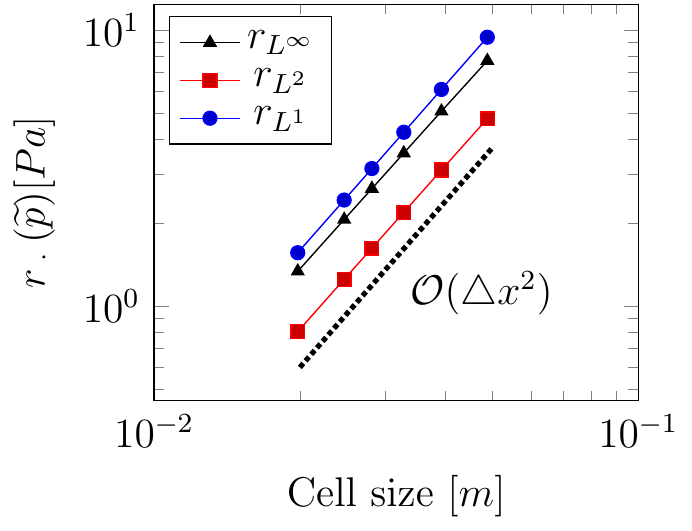}}
}
\caption{Error measurements for (a) velocity and (b) pressure of the transient three-dimensional  Example \ref{ex:tran3d}.}
\label{fig:transient_plots_3d}
\end{figure}

The solutions of the VANSE are chosen for time independent and time dependent cases, constructed similarly to ones of Blais et al.~\cite{BLAIS2015121} and Höcker et al.~\cite{hocker}. 
The here tested configurations are summarized as follows.
\begin{enumerate}
\item\label{ex:stat2d} Stationary two-dimensional example:
\begin{linenomath}\begin{align}
    \phi^{\flat} &= 0.5 + 0.4 \sin{(\pi x)} \sin{(\pi y)}, \\
    \overline{\bm{u}^{\flat}} &= 2 \begin{pmatrix} -(\sin{(\pi x)})^{2} \sin{(\pi y)} \cos{(\pi y)} \\ (\sin{(\pi y)})^{2} \sin{(\pi x)} \cos{(\pi x)} \end{pmatrix}, \\
    \tilde{p^{\flat}} &= \sin{(\pi x)} \sin{(\pi y)} .
\end{align}\end{linenomath}
\item\label{ex:stat3d} Stationary three-dimensional example:
\begin{linenomath}\begin{align}
    \phi^{\flat} &= 0.5 + 0.4 \sin{(\pi x)} \sin{(\pi y)} \sin{(\pi z)}, \\
    \overline{\bm{u}^{\flat}} &= \begin{pmatrix} (\sin{(\pi x)})^{2} \sin{(\pi y)} \cos{(\pi y)} \sin{(\pi z)} \cos{(\pi z)} \\ (\sin{(\pi y)})^{2} \sin{(\pi x)} \cos{(\pi x)} \sin{(\pi z)} \cos{(\pi z)} \\ -2 (\sin{(\pi z)})^{2} \sin{(\pi x)} \cos{(\pi x)} \sin{(\pi y)} \cos{(\pi y)} \end{pmatrix}, \\
    \tilde{p^{\flat}} &= \sin{(\pi x)} \sin{(\pi y)} \sin{(\pi z)} .
\end{align}\end{linenomath}
\item\label{ex:tran1d} Transient one-dimensional example:
\begin{linenomath}\begin{align}
    \phi^{\flat} &= 0.5 + 0.4 \sin{(\pi (x - 0.5t))}, \\
    \overline{\bm{u}^{\flat}} &= \begin{pmatrix} 0.5 + \frac{1}{\phi^{\flat}} \\ 0 \end{pmatrix}, \\
    \tilde{p^{\flat}} &= \sin{(\pi (x - 0.5t))} .
\end{align}\end{linenomath}
\item\label{ex:tran2d} Transient two-dimensional example:
\begin{linenomath}\begin{align}
    \phi^{\flat} &= 0.5 + 0.4 \sin{(\pi (x - 0.5t))} \sin{(\pi (y - 0.5t))}, \\
    \overline{\bm{u}^{\flat}} &= \begin{pmatrix} 0.5 + \frac{1}{\phi^{\flat}} \\ 0.5 + \frac{1}{\phi^{\flat}} \end{pmatrix}, \\
    \tilde{p^{\flat}} &= \sin{(\pi (x - 0.5t))} \sin{(\pi (y - 0.5t))} .
\end{align}\end{linenomath}
\item\label{ex:tran3d} Transient three-dimensional example:
\begin{linenomath}\begin{align}
    \phi^{\flat} &= 0.5 + 0.4 \sin{(\pi x)} \sin{(\pi y)} \sin{(\pi z)}, \\
    \overline{\bm{u}^{\flat}} &= \begin{pmatrix} 0.5 + \frac{1}{\phi^{\flat}} \\ 0.5 + \frac{1}{\phi^{\flat}} \\ 0.5 + \frac{1}{\phi^{\flat}} \end{pmatrix}, \\
    \tilde{p^{\flat}} &= \sin{(\pi x)} \sin{(\pi y)} \sin{(\pi z)} .
\end{align}\end{linenomath}
\end{enumerate}
The spatial simulation domain comprises $2\, m$ in each coordinate direction with periodic boundary conditions in every example. 
The fluid density is set to \(1\) $kg/m^{3}$ and the kinematic viscosity to \(0.1\) $m^{2}/s$. 
The relaxation time is held constant by all resolutions under diffusive scaling and is equal to \(0.53\) for the stationary and \(0.5075\) for the transient examples. 
Exemplary solutions are visualized in Figure \ref{fig:3D_stationary} for the stationary three-dimensional example \ref{ex:stat3d} and in Figure \ref{fig:3D_transient} for the transient three-dimensional example \ref{ex:tran3d}. 
The convergence plots for the examples in each error norm are shown in Figures \ref{fig:stationary_plots_2d}, \ref{fig:stationary_plots_3d}, \ref{fig:transient_plots_1d}, \ref{fig:transient_plots_2d}, and  \ref{fig:transient_plots_3d}, respectively. 

All examples are evaluated after the state stabilized and error norms remained asymptotically constant. 
It is to be noted that also the error norms reach a steady state after sufficiently long simulation time. 
This is due to the periodic boundary condition and constant maximal and minimal variable values that change only in position but not the amplitude.

In Figures \ref{fig:stationary_plots_2d}, \ref{fig:stationary_plots_3d}, \ref{fig:transient_plots_1d}, \ref{fig:transient_plots_2d}, and  \ref{fig:transient_plots_3d} we observe the same experimental convergence order of two in every norm type for the velocity as well as the pressure error. 
The absolute pressure deviation is not scaling by an increase of the target pressure values, so that the relative pressure can be made small enough. 
The results above clarify that the proposed LBM model for approximating VANSE converges with second order and thus is validly consistent in the present numerical tests.

\section{Conclusion}\label{sec:conc}
\noindent We establish a novel LBM for approximating the VANSE. 
The present LBM is formulated with an appropriate equilibrium distribution and pressure correction forcing term. 
The new moments of these equilibrium function and forcing terms, which take into account the local and temporal varying void fractions, are provided and justified. 
This unconventional point of view is based on considering streaming of the effective density from cell to cell. 
In particular, the population moments taken at one lattice cell include the effective density streamed from the neighbor cell in the chosen direction, so that a finite differences scheme is applicable. 

The numerical validation of the proposed LBM is performed on steady and transient examples, which are composed with MMS. 
Under the premise of diffusive scaling by refinement of the lattice resolution, the second order convergence of the fluid velocity and the pressure is approved.

Finally, the presented CE expansion formally validates the pressure correction forcing term via cancellation of moments with corresponding terms. 
Based on that, the expansion recovers the VANSE. 

In future studies the proposed LBM is to be extended to a full multiphase Eulerian model with phase interaction forces. 
Due to the intrinsic computing efficiency and optimal parallelizability of LBM, large eddy simulations \cite{smagorinsky} of complex entire reactor geometries with Eulerian multiphase LBM will become feasible. 
A second necessary extension of the model is the accounting for mass transfer between phases. 
Conclusively, the planned future research might render the multiphase LBM to an equal competitor of common FVM which is typically used in industrial solvers.

\section*{Author contribution}\noindent
\textbf{Fedor Bukreev}: 
Conceptualization,
Validation,
Formal analysis,
Investigation,
Resources,
Data Curation,
Writing - Original Draft;
\textbf{Stephan Simonis}: 
Methodology,
Validation, 
Formal analysis, 
Investigation,
Data Curation, 
Writing - Review \& Editing, 
Supervision; 
\textbf{Adrian Kummerländer}: 
Software,
Supervision;
\textbf{Julius Jeßberger}: 
Writing - Review \& Editing;
\textbf{Mathias J. Krause}: 
Software, 
Resources, 
Funding acquisition.

\section*{Acknowledgment}
\noindent
This work was performed on the HoreKa supercomputer funded by the Ministry of Science, Research and the Arts Baden-Württemberg and by the Federal Ministry of Education and Research. The current research is a part of the DFG project number 436212129 "Increase of efficiency in phosphate recovery by understanding the interaction of flow and loading processes with modeling and simulation".

\appendix
\section{Chapman\textendash Enskog analysis}
\label{sec:appendix}
\noindent
Below, we formally prove consistency of the above proposed LBE \eqref{eq:lbe} w.r.t. the targeted VANSE \eqref{eq:vansEquMass} and \eqref{eq:vansEquMom} up to higher order terms. 
The following CE expansion is based on the classical results for the diffusion limit towards the incompressible NSE as summarized for example in \cite{kruger} and references therein.

Let \(\epsilon > 0\) denote a label parameter for the Knudsen number \(\mathrm{Kn}\), and all other quantities be defined as above. 
We make the expansion ansatz
\begin{linenomath}\begin{align}\label{eq:pertAns}
\begin{cases}
    f_{i} = \sum\limits_{k=0}^{\infty} \epsilon^{k} f_{i}^{(k)} ,\\
    \partial_{t} = \epsilon \partial_{t}^{(1)} + \epsilon^{2} \partial_{t}^{(2)}, \\
    \bm{\nabla} = \epsilon \bm{\nabla}^{(1)}, \\
    \Omega_{i}^{F} = \epsilon \Omega_{i}^{F,(1)}.
\end{cases}
\end{align}\end{linenomath}
Taylor expanding the LBE \eqref{eq:lbe} yields
\begin{linenomath}\begin{align}\label{eq:lbeTaylored}
     \triangle t (\partial_{t} + \bm{\xi}_{i} \cdot \bm{\nabla})f_{i} + \frac{\triangle t^{2}}{2} (\partial_{t} + \bm{\xi}_{i} \cdot \bm{\nabla})^{2} f_{i} = & - \frac{\triangle t}{\tau} f_{i}^{\mathrm{neq}} + \triangle t \Omega_{i}^{F} \nonumber \\
     & + \mathcal{O}(\triangle t^3),
\end{align}\end{linenomath}
where \(f_{i}^\mathrm{neq} = f_{i} - f_{i}^{\mathrm{eq}}\). 
Subsequent to injecting \eqref{eq:pertAns}, the resulting version of \eqref{eq:lbeTaylored} can be separated into different \(\epsilon\)-orders, respectively

\noindent \(\mathcal{O}(\epsilon^{0})\):
\begin{linenomath}\begin{align}
    f_{i}^{(0)} = f_{i}^{eq}, 
\end{align}\end{linenomath}
\(\mathcal{O}(\epsilon^{1})\):
\begin{linenomath}\begin{align}
    (\partial_{t}^{(1)} + \bm{\xi}_{i} \cdot \bm{\nabla}^{(1)})f_{i}^{(0)} = - \frac{1}{\tau} f_{i}^{(1)} + \Omega_{i}^{F,(1)}, \label{eq:A1} 
\end{align}\end{linenomath}
\(\mathcal{O}(\epsilon^{2})\): 
\begin{linenomath}\begin{align}
    \partial_{t}^{(2)} & f_{i}^{(0)} + \Bigl(1 - \frac{\triangle t}{2 \tau} \Bigr) (\partial_{t}^{(1)} + \bm{\xi}_{i} \cdot \bm{\nabla}^{(1)}) \Bigl( f_{i}^{(1)} + \frac{\triangle t}{2} \Omega_{i}^{F,(1)} \Bigr) = - \frac{1}{\tau} f_{i}^{(2)}. \label{eq:A2}
\end{align}\end{linenomath}
Here and in the following, we neglect derivative terms of order \(k\geq 3\) due to the smallness argument \(\triangle t^{k} (\partial_{t} + \bm{\epsilon}_{i}\cdot \bm{\nabla})^{k} f_{i} \sim \mathcal{O} (\mathrm{Kn}^{k})\), as explained in detail in \cite{kruger}. 

The pressure forcing term moments, taken in one lattice cell, are
\begin{linenomath}\begin{align}
    M_0^{F,(1)} &= \sum_i \Omega_i^{F_{\mathrm{PC}}} = 0, \label{eq:M0F}\\
    M_{1 \alpha}^{F,(1)} &= \sum_i \xi_{i \alpha} \Omega_i^{F_{\mathrm{PC}}} = \Bigl(1 - \frac{\triangle t}{2 \tau} \Bigr) F_{\mathrm{PC}, \alpha} = \Bigl(1 - \frac{\triangle t}{2 \tau} \Bigr) p \bm{\nabla} \phi, \\
     M_{2 \alpha \beta}^{F,(1)} &= \sum_i \xi_{i \alpha} \xi_{i \beta} \Omega_i^{F_{\mathrm{PC}}} = \Bigl(1 - \frac{\triangle t}{2 \tau} \Bigr) (F_{\mathrm{PC}, \alpha} u_{\beta} + u_{\alpha} F_{\mathrm{PC}, \beta}). \label{eq:M2F}
\end{align}\end{linenomath}

Note that
\begin{linenomath}\begin{align}
    M_{0}^{(1)} &= - \frac{\triangle t}{2\left( 1 - \frac{\triangle t}{2 \tau}\right)} M_{0}^{F,(1)} = 0,  \label{eq:ceExpBaseI}\\
    M_{1\alpha}^{(1)} &= - \frac{\triangle t}{2\left( 1 - \frac{\triangle t}{2 \tau}\right)} M_{1\alpha}^{F,(1)} = -\frac{\triangle t}{2} F_{\mathrm{PC}, \alpha}, \label{eq:ceExpBaseII}
\end{align}\end{linenomath}
and for \(k\geq 2\) is assumed
\begin{linenomath}\begin{align} \label{eq:ceExpBaseIII}
    M_{0}^{(k)} = \sum\limits_{i} f_{i}^{(k)} = 0 =  \sum\limits_{i} \xi_{i\alpha} f_{i}^{(k)} = M_{1\alpha}^{(k)}.
\end{align}\end{linenomath}
Taking the zeroth, first and second order moments of \eqref{eq:A1} and the zeroth and first ones of \eqref{eq:A2}, and substituting the notation of (\ref{eq:M0}\textendash \ref{eq:M3}) and (\ref{eq:M0F}\textendash \ref{eq:M2F}), we obtain respectively 

\noindent $\mathcal{O}(\epsilon^{1})$:
\begin{linenomath}\begin{align}
    \partial_{t}^{(1)} M_{0}^{(0)} + \bm{\nabla}^{(1)} M_{1\alpha}^{(0)} &= 0, \label{eq:mom0Eps1}\\
    \partial_{t}^{(1)} M_{1\alpha}^{(0)} + \bm{\nabla}^{(1)} M_{2\alpha\beta}^{(0)} &= -\frac{1}{\tau}M_{1\alpha}^{(1)} + M_{1\alpha}^{F,(1)}, \label{eq:mom1Eps1} \\
    \partial_{t}^{(1)} M_{2\alpha\beta}^{(0)} + \bm{\nabla}^{(1)} M_{3\alpha\beta\gamma}^{(0)} &= -\frac{1}{\tau} M_{2\alpha\beta}^{(1)} + M_{2\alpha\beta}^{F,(1)} \label{eq:mom2Eps1}
\end{align}\end{linenomath}
$\mathcal{O}(\epsilon^{2})$:
\begin{linenomath}\begin{align}
    \partial_{t}^{(2)} M_{0}^{(0)} &= 0, \label{eq:mom0Eps2} \\
    \partial_{t}^{(2)} M_{1\alpha}^{(0)} + \left(1 - \frac{\triangle t}{2 \tau} \right) \bm{\nabla}^{(1)} M_{2\alpha\beta}^{(1)} &= - \frac{\triangle t}{2}  \bm{\nabla}^{(1)} \left(1 - \frac{\triangle t}{2 \tau} \right) M_{2\alpha\beta}^{F,(1)}. \label{eq:mom1Eps2}
\end{align}\end{linenomath}
Thus, the recombination \eqref{eq:mom0Eps1} \(+\) \eqref{eq:mom0Eps2} and \eqref{eq:mom1Eps1} \(+\) \eqref{eq:mom1Eps2} yields
\begin{linenomath}\begin{align}\label{eq:mom0EpsSum}
    \partial_{t} M_{0}^{(0)} + \bm{\nabla}^{(1)} M_{1\alpha}^{(0)} 
    = 0
\end{align}\end{linenomath}
and
\begin{linenomath}\begin{align}\label{eq:mom1EpsSum}
       \partial_{t} & M_{1\alpha}^{(0)} + \bm{\nabla}^{(1)} M_{2\alpha\beta}^{(0)} + \left(1 - \frac{\triangle t}{2 \tau} \right) \bm{\nabla}^{(1)} M_{2\alpha\beta}^{(1)} \nonumber \\
       &= -\frac{1}{\tau}M_{1\alpha}^{(1)} + M_{1\alpha}^{F,(1)} - \frac{\triangle t}{2} \left(1 - \frac{\triangle t}{2 \tau} \right) \bm{\nabla}^{(1)} M_{2\alpha\beta}^{F,(1)} ,
\end{align}\end{linenomath}
respectively. 
Under the diffusion limit assumption $\triangle t \sim \triangle x^{2} \to 0$ when refining the spatial mesh \(\triangle x \to 0\), we rewrite \eqref{eq:mom0EpsSum} and \eqref{eq:mom1EpsSum} as
\begin{linenomath}\begin{align}
    \partial_{t} M_{0}^{(0)} + \bm{\nabla}^{(1)} M_{1\alpha}^{(0)} &= 0, \label{eq:A3} \\
    \partial_{t} M_{1\alpha}^{(0)} + \bm{\nabla}^{(1)} M_{2\alpha\beta}^{(0)} &= -\frac{1}{\tau}M_{1\alpha}^{(1)} + M_{1\alpha}^{F,(1)}  - \left(1 - \frac{\triangle t}{2 \tau} \right) \bm{\nabla}^{(1)} M_{2\alpha\beta}^{(1)}. \label{eq:A4}
\end{align}\end{linenomath}
We substitute the moment notation to recover 
\begin{linenomath}\begin{align}
   & \partial_t \left( \tilde{\rho^{\flat}} \Bigl(\phi^{\flat} (\bm{x},t) + \varpi_{i \neq 0}(N) \bm{\nabla}^{2} \phi^{\flat} \Bigr) \right) \nonumber \\
   &  + \bm{\nabla} \cdot \left(\tilde{\rho^{\flat}} \Bigl(\phi^{\flat} (\bm{x},t) + \varpi_{i \neq 0}(N) \bm{\nabla}^{2} \phi^{\flat} \Bigr) \overline{\bm{u}^{\flat}} \right) && = 0, \\
   & \partial_t \left( \tilde{\rho^{\flat}} \Bigl(\phi^{\flat} (\bm{x},t) + \varpi_{i \neq 0}(N) \bm{\nabla}^{2} \phi^{\flat} \Bigr) \overline{\bm{u}^{\flat}} \right) \nonumber \\
   & + \bm{\nabla} \cdot \left( \tilde{\rho^{\flat}} \Bigl(\phi^{\flat} (\bm{x},t) + \varpi_{i \neq 0}(N) \bm{\nabla}^{2} \phi^{\flat} \Bigr) \overline{\bm{u}^{\flat}} \overline{\bm{u}^{\flat}} \right) \nonumber \\ 
   & + \bm{\nabla} \left( \tilde{\rho^{\flat}} c_s^{2} \Bigl(\phi^{\flat} (\bm{x},t) + \varpi_{i \neq 0}(N) \bm{\nabla}^{2} \phi^{\flat} \Bigr) \right) && = -\Bigl(1 - \frac{\triangle t}{2 \tau} \Bigr) \bm{\nabla} M_{2\alpha\beta}^{(1)} \nonumber \\
   &  && \hphantom{=~} + \tilde{\rho^{\flat}} c_{s}^{2} \bm{\nabla} \phi^{\flat}.
\end{align}\end{linenomath}
Via reordering \eqref{eq:mom2Eps1} and deletion of higher order terms, we unfold 

\begin{linenomath}\begin{align}
    M_{2 \alpha \beta}^{(1)} = c_{s}^{2} \phi^{\flat} \tilde{\rho^{\flat}} (\bm{\nabla} \overline{\bm{u}^{\flat}} + \overline{\bm{u}^{\flat}} \bm{\nabla}).
\end{align}\end{linenomath}
After deletion of vanishing terms of the 3rd order and with $\triangle t$ and insertion of the stress tensor $\bm{\nabla} M_{2\alpha\beta}^{(1)}$ the VANSE are recovered up to higher order terms 
\begin{linenomath}\begin{align}
    \partial_t (\phi^{\flat} \tilde{\rho^{\flat}})+ \bm{\nabla} \cdot (\phi^{\flat} \tilde{\rho^{\flat}} \overline{\bm{u}^{\flat}}) &= 0, \\
    \partial_t (\phi^{\flat} \tilde{\rho^{\flat}} \overline{\bm{u}^{\flat}})+ \bm{\nabla} \cdot (\phi^{\flat} \tilde{\rho^{\flat}} \overline{\bm{u}^{\flat}} \overline{\bm{u}^{\flat}})+ \phi^{\flat} \bm{\nabla} \tilde{p} &= \nu \bm{\nabla} \cdot (\phi^{\flat} \tilde{\rho^{\flat}} (\bm{\nabla} \overline{\bm{u}^{\flat}} + \overline{\bm{u}^{\flat}} \bm{\nabla}))
\end{align}\end{linenomath}
where the viscosity is regained as
\begin{linenomath}\begin{align}
    \nu = \Bigl(\tau - \frac{\triangle t}{2} \Bigr) c_{s}^{2}, 
\end{align}\end{linenomath}
from comparison to \eqref{eq:vansEquMom}.


\begin{thebibliography}{10}
\expandafter\ifx\csname url\endcsname\relax
  \def\url#1{\texttt{#1}}\fi
\expandafter\ifx\csname urlprefix\endcsname\relax\def\urlprefix{URL }\fi
\expandafter\ifx\csname href\endcsname\relax
  \def\href#1#2{#2} \def\path#1{#1}\fi

\bibitem{reschetilowski}
W.~Reschetilowski, Handbuch Chemische Reaktoren Grundlagen und Anwendungen der
  Chemischen Reaktionstechnik: Grundlagen und Anwendungen der Chemischen
  Reaktionstechnik, Springer Spektrum, Berlin Heidelberg, 2020.
\newblock \href {https://doi.org/10.1007/978-3-662-56444-8}
  {\path{doi:10.1007/978-3-662-56444-8}}.

\bibitem{sattler}
K.~Sattler, Thermische Trennverfahren, John Wiley {\&} Sons, Ltd., 2001.
\newblock \href {https://doi.org/10.1002/3527603328.ch1a}
  {\path{doi:10.1002/3527603328.ch1a}}.

\bibitem{bohnet}
M.~Bohnet, Mechanische Verfahrenstechnik, John Wiley {\&} Sons, Ltd, 2003.
\newblock \href {https://doi.org/10.1002/9783527663569.fmatter}
  {\path{doi:10.1002/9783527663569.fmatter}}.

\bibitem{hiltunen}
K.~Hiltunen, A.~Jäsberg, S.~Kallio, H.~Karema, M.~Kataja, A.~Koponen,
  M.~Manninen, V.~Taivassalo, {Multiphase flow dynamics: Theory and numerics},
  VTT Publications 722, VTT Technical Research Centre of Finland, 2009.

\bibitem{GIDASPOW1994337}
D.~Gidaspow, Multiphase Flow and Fluidization, Academic Press, San Diego, 1994.
\newblock \href {https://doi.org/10.1016/C2009-0-21244-X}
  {\path{doi:10.1016/C2009-0-21244-X}}.

\bibitem{zhu}
T.~Zhu, \href{https://mediatum.ub.tum.de/doc/1279870/1279870.pdf}{{Unsteady
  porous\textendash media flow}}, Ph.D. thesis, Technische Universität
  München (2016).
\newline\urlprefix\url{https://mediatum.ub.tum.de/doc/1279870/1279870.pdf}

\bibitem{PEPIOT2012104}
P.~Pepiot, O.~Desjardins, {Numerical analysis of the dynamics of two- and
  three-dimensional fluidized bed reactors using an Euler–Lagrange approach},
  Powder Technology 220 (2012) 104--121.
\newblock \href {https://doi.org/10.1016/j.powtec.2011.09.021}
  {\path{doi:10.1016/j.powtec.2011.09.021}}.

\bibitem{moukalled}
F.~Moukalled, L.~Mangani, M.~Darwish, The Finite Volume Method in Computational
  Fluid Dynamics. An Advanced Introduction with OpenFOAM and Matlab, Springer
  International Publishing Switzerland, 2016.
\newblock \href {https://doi.org/10.1007/978-3-319-16874-6}
  {\path{doi:10.1007/978-3-319-16874-6}}.

\bibitem{kruger}
T.~Krüger, H.~Kusumaatmaja, A.~Kuzmin, O.~Shardt, G.~Silva, E.~M. Viggen, The
  Lattice Boltzmann Method - Principles and Practice, Springer Cham, 2016.
\newblock \href {https://doi.org/10.1007/978-3-319-44649-3}
  {\path{doi:10.1007/978-3-319-44649-3}}.

\bibitem{hanel2004molekulare}
D.~Hänel, {Molekulare Gasdynamik : Einführung in die kinetische Theorie der
  Gase und Lattice-Boltzmann-Methoden}, Springer, Berlin Heidelberg, 2004.
\newblock \href {https://doi.org/10.1007/3-540-35047-0}
  {\path{doi:10.1007/3-540-35047-0}}.

\bibitem{li}
J.~Li, \href{https://arxiv.org/abs/1512.02599}{{Appendix: Chapman-Enskog
  Expansion in the Lattice Boltzmann Method}} (2015).
\newblock \href {https://doi.org/10.48550/ARXIV.1512.02599}
  {\path{doi:10.48550/ARXIV.1512.02599}}.
\newline\urlprefix\url{https://arxiv.org/abs/1512.02599}

\bibitem{simonis2022limit}
S.~Simonis, M.~J. Krause, \href{https://arxiv.org/abs/2208.06867}{{Limit consistency of lattice Boltzmann equations}} (2022).
\newblock \href {https://doi.org/10.48550/arXiv.2208.06867}
  {\path{doi:10.48550/arXiv.2208.06867}}.
\newline\urlprefix\url{https://arxiv.org/abs/2208.06867}

\bibitem{simonis2022temporal}
S.~Simonis, D.~Oberle, M.~Gaedtke, P.~Jenny, M.~J. Krause, {Temporal large eddy
  simulation with lattice Boltzmann methods}, Journal of Computational Physics
  454 (2022) 110991.
\newblock \href {https://doi.org/10.1016/j.jcp.2022.110991}
  {\path{doi:10.1016/j.jcp.2022.110991}}.

\bibitem{simonis2021linear}
S.~Simonis, M.~Haussmann, L.~Kronberg, W.~D{\"o}rfler, M.~J. Krause, {Linear
  and brute force stability of orthogonal moment multiple-relaxation-time
  lattice Boltzmann methods applied to homogeneous isotropic turbulence},
  Philosophical Transactions of the Royal Society A 379~(2208) (2021) 20200405.
\newblock \href {https://doi.org/10.1098/rsta.2020.0405}
  {\path{doi:10.1098/rsta.2020.0405}}.

\bibitem{haussmann2019direct}
M.~Haussmann, S.~Simonis, H.~Nirschl, M.~J. Krause, {Direct numerical
  simulation of decaying homogeneous isotropic turbulence—numerical
  experiments on stability, consistency and accuracy of distinct lattice
  Boltzmann methods}, International Journal of Modern Physics C 30~(09) (2019)
  1950074.
\newblock \href {https://doi.org/10.1142/S0129183119500748}
  {\path{doi:10.1142/S0129183119500748}}.

\bibitem{simonis2020relaxation}
S.~Simonis, M.~Frank, M.~J. Krause, {On relaxation systems and their relation
  to discrete velocity Boltzmann models for scalar advection--diffusion
  equations}, Philosophical Transactions of the Royal Society A 378~(2175)
  (2020) 20190400.
\newblock \href {https://doi.org/10.1098/rsta.2019.0400}
  {\path{doi:10.1098/rsta.2019.0400}}.

\bibitem{dapelo2021lattice-boltzmann}
D.~Dapelo, S.~Simonis, M.~J. Krause, J.~Bridgeman, {Lattice-Boltzmann coupled
  models for advection–diffusion flow on a wide range of Péclet numbers},
  Journal of Computational Science 51 (2021) 101363.
\newblock \href {https://doi.org/10.1016/j.jocs.2021.101363}
  {\path{doi:10.1016/j.jocs.2021.101363}}.

\bibitem{mink2021comprehensive}
A.~Mink, K.~Schediwy, C.~Posten, H.~Nirschl, S.~Simonis, M.~J. Krause,
  \href{https://arxiv.org/abs/2107.12210}{{Comprehensive computational model
  for coupled fluid flow, mass transfer and light supply in tubular
  photobioreactors equipped with glass sponges}} (2021).
\newblock \href {https://doi.org/10.48550/ARXIV.2107.12210}
  {\path{doi:10.48550/ARXIV.2107.12210}}.
\newline\urlprefix\url{https://arxiv.org/abs/2107.12210}

\bibitem{simonis2022forschungsnahe}
S.~Simonis, M.~J. Krause, {Forschungsnahe Lehre unter Pandemiebedingungen},
  Mitteilungen der Deutschen Mathematiker-Vereinigung 30~(1) (2022) 43--45.
\newblock \href {https://doi.org/10.1515/dmvm-2022-0015}
  {\path{doi:10.1515/dmvm-2022-0015}}.

\bibitem{haussmann2021fluid-structure}
M.~Haussmann, P.~Reinshaus, S.~Simonis, H.~Nirschl, M.~J. Krause,
  {Fluid{\textendash}Structure Interaction Simulation of a Coriolis Mass
  Flowmeter Using a Lattice Boltzmann Method}, Fluids 6~(4) (2021).
\newblock \href {https://doi.org/10.3390/fluids6040167}
  {\path{doi:10.3390/fluids6040167}}.

\bibitem{krause}
M.~Krause, A.~Kummerländer, S.~Avis, H.~Kusumaatmaja, D.~Dapelo, F.~Klemens,
  M.~Gaedtke, N.~Hafen, A.~Mink, R.~Trunk, J.~Marquardt, M.-L. Maier,
  M.~Haussmann, S.~Simonis, {OpenLB{\textemdash}Open source lattice Boltzmann
  code}, Computers {\&} Mathematics with Applications 81 (2020) 258--288.
\newblock \href {https://doi.org/10.1016/j.camwa.2020.04.033}
  {\path{doi:10.1016/j.camwa.2020.04.033}}.

\bibitem{kummerlander2022olb15}
A.~Kummerl\"ander, S.~Avis, H.~Kusumaatmaja, F.~Bukreev, D.~Dapelo,
  S.~Gro\ss{}mann, N.~Hafen, C.~Holeksa, A.~Husfeldt, J.~Je\ss{}berger,
  L.~Kronberg, J.~Marquardt, J.~M\"odl, J.~Nguyen, T.~Pertzel, S.~Simonis,
  L.~Springmann, N.~Suntoyo, D.~Teutscher, M.~Zhong, M.~Krause,
  \href{https://doi.org/10.5281/zenodo.6469606}{{OpenLB Release 1.5: Open
  Source Lattice Boltzmann Code}} (Nov. 2022).
\newblock \href {https://doi.org/10.5281/zenodo.6469606}
  {\path{doi:10.5281/zenodo.6469606}}.
\newline\urlprefix\url{https://doi.org/10.5281/zenodo.6469606}

\bibitem{tuprints}
M.~Haussmann, F.~Ries, J.~B. Jeppener-Haltenhoff, Y.~Li, M.~Schmidt, C.~Welch,
  L.~Illmann, B.~B{\"o}hm, H.~Nirschl, M.~J. Krause, A.~Sadiki, {Evaluation of
  a Near-Wall-Modeled Large Eddy Lattice Boltzmann Method for the Analysis of
  Complex Flows Relevant to IC Engines}, Computation 8~(2) (2020).
\newblock \href {https://doi.org/10.25534/tuprints-00013372}
  {\path{doi:10.25534/tuprints-00013372}}.

\bibitem{guo}
Z.~Guo, T.~Zhao, {Lattice Boltzmann model for incompressible flows through
  porous media}, Physical Review E 66 (2002) 036304.
\newblock \href {https://doi.org/10.1103/PhysRevE.66.036304}
  {\path{doi:10.1103/PhysRevE.66.036304}}.

\bibitem{blais}
B.~Blais, J.-M. Tucny, D.~Vidal, F.~Bertrand, {A conservative lattice Boltzmann
  model for the volume-averaged Navier–Stokes equations based on a novel
  collision operator}, Journal of Computational Physics 294 (2015) 258--273.
\newblock \href {https://doi.org/10.1016/j.jcp.2015.03.036}
  {\path{doi:10.1016/j.jcp.2015.03.036}}.

\bibitem{hocker}
S.~Höcker, R.~Trunk, W.~Dörfler, M.~Krause, {Towards the Simulations of
  Inertial Dense Particulate Flows with a Volume-Averaged Lattice Boltzmann
  Method}, Computers {\&} Fluids 166 (2018) 152--162.
\newblock \href {https://doi.org/10.1016/j.compfluid.2018.02.011}
  {\path{doi:10.1016/j.compfluid.2018.02.011}}.

\bibitem{Maier2021_1000132643}
M.-L. Maier, {Coupled lattice Boltzmann and discrete element method for
  reactive particle fluid flows with applications in process engineering},
  Ph.D. thesis, Karlsruher Institut für Technologie (KIT) (2021).
\newblock \href {https://doi.org/10.5445/IR/1000132643}
  {\path{doi:10.5445/IR/1000132643}}.

\bibitem{zhang}
J.~Zhang, L.~Wang, J.~Ouyang, {Lattice Boltzmann Model for The Volume-Averaged
  Navier-Stokes Equations}, EPL (Europhysics Letters) 107 (2014) 20001.
\newblock \href {https://doi.org/10.1209/0295-5075/107/20001}
  {\path{doi:10.1209/0295-5075/107/20001}}.

\bibitem{bgk}
P.~L. Bhatnagar, E.~P. Gross, M.~Krook, {A Model for Collision Processes in
  Gases. I. Small Amplitude Processes in Charged and Neutral One-Component
  Systems}, Physical Review E 94 (1954) 511--525.
\newblock \href {https://doi.org/10.1103/PhysRev.94.511}
  {\path{doi:10.1103/PhysRev.94.511}}.

\bibitem{guozhaoli}
Z.~Guo, C.~Zheng, B.~Shi, {Discrete lattice effects on the forcing term in the
  lattice Boltzmann method}, Physical Review E 65 (2002) 046308.
\newblock \href {https://doi.org/10.1103/PhysRevE.65.046308}
  {\path{doi:10.1103/PhysRevE.65.046308}}.

\bibitem{roache}
P.~Roache, {Code Verification by the Method of Manufactured Solutions}, Journal
  of Fluids Engineering 124 (2002) 4.
\newblock \href {https://doi.org/10.1115/1.1436090}
  {\path{doi:10.1115/1.1436090}}.

\bibitem{oberkampf}
W.~Oberkampf, C.~Roy, Verification and Validation in Scientific Computing,
  Cambridge University Press, 2010.
\newblock \href {https://doi.org/10.1017/CBO9780511760396}
  {\path{doi:10.1017/CBO9780511760396}}.

\bibitem{BLAIS2015121}
B.~Blais, F.~Bertrand, {On the use of the method of manufactured solutions for
  the verification of CFD codes for the volume-averaged Navier–Stokes
  equations}, Computers {\&} Fluids 114 (2015) 121--129.
\newblock \href {https://doi.org/10.1016/j.compfluid.2015.03.002}
  {\path{doi:10.1016/j.compfluid.2015.03.002}}.

\bibitem{smagorinsky}
J.~Smagorinsky, {General circulation experiments with the primitive equations:
  I. The basic experiment}, Monthly Weather Review 91~(3) (1963) 99--164.
\newblock \href
  {https://doi.org/10.1175/1520-0493(1963)091<0099:GCEWTP>2.3.CO;2}
  {\path{doi:10.1175/1520-0493(1963)091<0099:GCEWTP>2.3.CO;2}}.

\end{thebibliography}
\end{document}